# A Variable Reduction Method for Large-Scale Security Constrained Unit Commitment

Xuan Li, Qiaozhu Zhai, *Member, IEEE*, Jingxuan Zhou, Xiaohong Guan, *Fellow, IEEE*

*Abstract*—Efficient methods for large-scale security constrained unit commitment (SCUC) problems have long been an important research topic and a challenge especially in market clearing computation. For large-scale SCUC, the Lagrangian relaxation methods (LR) and the mixed integer programming methods (MIP) are most widely adopted. However, LR usually suffers from slow convergence; and the computational burden of MIP is heavy when the binary variable number is large. In this paper, a variable reduction method is proposed: First, the time-coupled constraints in the original SCUC problem are relaxed and many single-period SCUC problems (s-UC) are obtained. Second, LR is used to solve the s-UCs. Different from traditional LR with iterative subgradient method, it is found that the optimal multipliers and the approximate UC solutions of s-UCs can be obtained by solving linear programs. Third, a criterion for choosing and fixing the UC variables in the SCUC problem is established, hence the number of binary variables is reduced. Last, the SCUC with reduced binary variables is solved by MIP solver to obtain the final UC solution. The proposed method is tested on the IEEE 118-bus system and a 6484-bus system. The results show the method is very efficient and effective.

## Nomenclature

### A. Problems

**f-UC**     Full-scale SCUC problem, as (1)-(7).

**s-UC**     Single-period SCUC, as (8)-(13).

**d-UC**     Dual problem of **s-UC**, as (18)-(19).

**d-UC-LP** linear program formulation for **d-UC**, as (38)-(42).

**i-UC**     Individual-unit subproblem, as (26)-(28).

### B. Indices and Sets

$i$     Index for units, with the number of $I$.

$k$     Index for loads.

$l$     Index for transmission capacity/security constraints, with the number of $L$.

$t$     Index for time periods, with the number of $T$.

$X$     Set of allowable decisions of binary variables.

This work is supported in part by National Key R&D Project (2016YFB0901900), National Natural Science Foundation (61773309, 61773308, 61473218) of China and by Open Fund of State Key Laboratory of Operation and Control of Renewable Energy & Storage Systems.

X. Li, Q. Zhai, J. Zhou, and X. Guan are with Systems Engineering Institute, MOEKLINNS Lab, Xi'an Jiaotong University, Xi'an 710049, China. (Correspondence author: qzzhai@sei.xjtu.edu.cn).

### C. Constants

$a_i, b_i$     Coefficients for the linear fuel cost function of unit $i$.

$d_{k,t}$     Load level of load $k$ at period $t$ (MW).

$\Delta_i^+, \Delta_i^-$     Ramp up/down limits of unit $i$ (MW).

$F_l$     Transmission limit in transmission capacity/security constraint $l$ (MW).

$\Gamma$     Power transfer distribution factor (PTDF).

$\underline{P}_i, \overline{P}_i$     Min/max generation capacity of unit $i$ (MW).

### D. Variables

$p_{i,t}$     Dispatch decision of unit $i$ at period $t$ (MW).

$z_{i,t}$     Unit commitment decision of unit $i$ at period $t$ (binary).

$\tilde{z}$     Approximate solution of **s-UC** obtained by solving **d-UC-LP**.

$z'$     Optimal solution of **s-UC**.

$\lambda_{0,t}$     Dual variable for power balance constraint at period $t$.

$\lambda_{l,t}^+, \lambda_{l,t}^-$     Dual variables for transmission capacity/security constraint $l$ at period $t$.

$\boldsymbol{\lambda}_t$     Vector of all dual variables at period $t$.

### E. Functions

$C(\cdot)$     Function of the fuel cost ($).

$S(\cdot)$     Function of the start-up cost ($).

## I. Introduction

Security constrained unit commitment problem (SCUC) is a fundamental problem for power system's planning, scheduling and operations [1,2]. SCUC is usually modeled as a mixed-integer optimization problem, which aims to determine a set of unit commitment (UC) and economic dispatch (ED) decisions with minimum total costs, while satisfying many kinds of system-wide and individual-unit constraints [2]. With the rapid development of the power systems, the problem scale and complexity of SCUC also increase significantly; more complex transmission network, more thermal units, variant energy sources and uncertainties must be considered and modeled in the problems (e.g., 45,000 buses, 1,400 generation resources and 2,446 pricing nodes are included in the MISO's model [3]). In this case, solving large-scale SCUC problem can be very challenging [6]. In addition, as an essential part of the day-ahead market clearing process, the solution of the large-scale SCUC must be efficient enough to meet the market clearing time requirements. (e.g., 3~4 hours in MISO, ISO New England and PJM [3-5]). Therefore, it is very important to quickly obtain a feasible solution (even the solution quality is not good enough) in case the solvers fail to return a solution in the re-quired time limit.

For large-scale SCUC problems, Lagrangian relaxation-based methods (LR) like [6-11] and general mixed integer programming-based methods (MIP) like [12-18] are most widely adopted [3,19].

The general solution process of the traditional LR-based methods is as follows: first the system-wide constraints, such as power balance and transmission capacity/security constraints are relaxed by introducing Lagrangian multipliers, and the dual problem is formed. Then an iterative framework is established to solve the dual problem. In each iteration, with given multipliers, the dual function (objective function of the dual problem) is obtained by solving a number of individual-unit subproblems. Afterward, the multipliers are updated using the subgradient-based method (or some other methods). Then, the subproblems are resolved with the updated multipliers, and so forth until convergence. The salient advantage of the LR based approach is that the computational complexity of solving the dual problem is linearly related to the unit number [20]. However, the convergence of such algorithm usually suffers from zigzagging, and is sometimes very slow [7]. Moreover, at convergence, the feasibility of the solution is not guaranteed. Methods like [8,20] must be adopted to construct a feasible solution.

To improve the performance of LR, [7] develops the surrogate LR method and proves its convergence to the optimal multipliers. This work is further applied to large-scale SCUC with combined cycle units in [6]. [9,10] propose enhanced adaptive LR methods with heuristics. [11] uses the augmented LR to solve the unit commitment problem and compares two decomposition approaches.

MIP-based methods are prevalent in recent years due to the development of the commercial MIP solver like CPLEX and Gurobi [14]. Compared with LR, solving SCUC by using MIP-based methods is more reliable and straightforward, and the problem formulation can be more complicated [3]. However, when large-scale SCUC problems are solved, the computational complexity of the MIP methods can be a serious issue.

To reduce the computational burden of MIP methods, many works focus on efficient modeling skills: [16] proposes a locally ideal formulation for piecewise linear cost functions. [12,13,15] propose computationally efficient formulations with one, two and three binary UC variables for each unit at each period, respectively. [17] tightens the start-up and shut-down ramping constraints in the UC problem. Large numbers of security constraints are seen as one of the main challenges when solving the SCUC problem. Some works focus on reducing the constraint number in the solution procedure: the method proposed in [21] can quickly identify most of the redundant security constraints with an analytical sufficient condition. In [3], two iterative methods are reported that only include binding or near binding constraints. [3] also reports a binary reduction approach that can utilize the existing UC solution and price information. [18] reduces the overall problem scale by decomposing the multi-period problems into single-period ones,

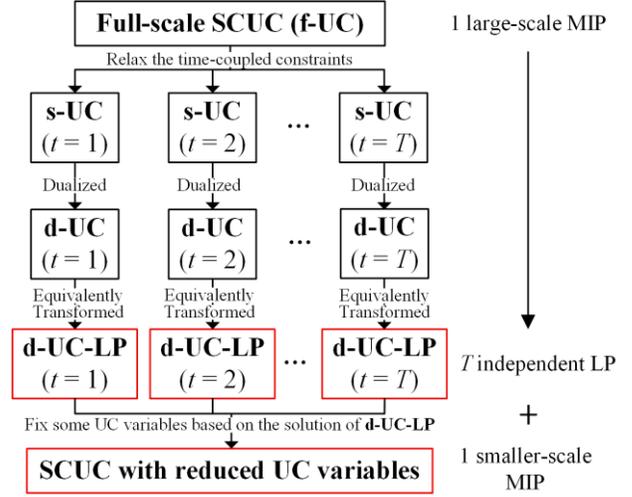

Fig. 1 Framework for solving large-scale SCUC problem.

and form the feasible solution within an iterative framework.

The MIP solvers are usually based on branch-and-cut procedure (BAC), the overall computational burden is closely related to 1) the number of nodes explored in the BAC tree, and 2) the time spent on exploring each node [16]. With large numbers of UC variables, the size of BAC tree can be very large, and the problem solved at each node can also be computationally intensive.

In this paper, a new solution framework for large-scale SCUC problem is proposed, the main idea is depicted in Fig. 1, and summarized as follows:

1) By relaxing the time-coupled constraints, the full-scale multi-period SCUC problem (f-UC) is transformed into $T$ single-period SCUC problems (s-UC).

2) For each period $t$, instead of directly applying the MIP solver to solve the s-UC problem. We dualize s-UC and obtain its Lagrangian dual problem d-UC.

3) The d-UC problems are then equivalently transformed into linear programs (LP) d-UC-LP. By solving d-UC-LP, the optimal multipliers and the approximate UC solutions of s-UC can be directly obtained. This is completely different from the conventional iterative subgradient-based method for the Lagrangian dual problem.

4) Based on the solution of d-UC-LP, a variable reduction criterion is established, by which some binary UC variables are fixed to the values of the approximate solution of s-UC. Then the f-UC with much reduced binary variables are solved by MIP solver to obtain the final solution.

The rest of this paper is organized as follows: section II presents the basic SCUC formulation. Section III delineates the variable reductions methods. Section IV presents the numerical results. Section V concludes the paper.

## II. PROBLEM FORMULATION

The SCUC problem is formulated in this section. The problem is denoted as f-UC (abbreviation for "full-scale" SCUC).

**(f-UC):**

$$\min_{z,p} \sum_i S_i(z_{i,1},...,z_{i,T}) + \sum_{i,t} C_i(p_{i,t},z_{i,t}) \qquad (1)$$

$$C_i(p_{i,t},z_{i,t}) = \begin{cases} 0; & \text{if } z_{i,t}=0 \\ a_i p_{i,t}+b_i; & \text{if } z_{i,t}=1 \end{cases} \qquad (2)$$

$$\sum_i p_{i,t} - \sum_k d_{k,t} = 0; \forall t \qquad (3)$$

$$-F_l \le \sum_i \Gamma_{l,i}^U p_{i,t} - \sum_k \Gamma_{l,k}^D d_{k,t} \le F_l; \forall l, t \qquad (4)$$

$$z_{i,t}\underline{P} \le p_{i,t} \le z_{i,t}\overline{P}_i; \forall i, t \qquad (5)$$

$$\Delta_i^-(z_{i,t},z_{i,t-1}) \le p_{i,t} - p_{i,t-1} \le \Delta_i^+(z_{i,t},z_{i,t-1}); \forall i, t \qquad (6)$$

$$z \in X, z \in \{0,1\}^{I \times T} \qquad (7)$$

The decision variables are the UC decision $z \in \{0,1\}^{I \times T}$ and economic dispatch (ED) decision $p \in \mathbb{R}^{I \times T}$.

The objective (1) is to minimize the total start-up cost and fuel cost over the scheduling horizon. For better presentation, the fuel cost function (2) is assumed linear in section II and III. Main results in this paper are still valid for piecewise-linear fuel cost functions, and some details are given in the Appendix.

The constraints include (3)-(7): (3) is the power balance constraint. (4) is the general form of transmission capacity/security constraint based on DC flow model. Suppose the number of transmission lines is $N$, then the number of transmission capacity constraints is $N$, and the number of "$N$-1" transmission security constraints is no more than $N(N$-1). $\Gamma$ in (4) is the power transfer distribution factor (PTDF) that represents the sensitivity of power flow on the transmission line to the nodal power injection. (5) is the generation capacity constraint. (6) represents the ramp up/down constraints, the detailed formulations of the ramp limits $\Delta^+$, $\Delta^-$ can be found in [16]. (7) represents the constraints only for UC decisions, including minimum up/down constraints, must on/off constraints, and so on. According to [22] (the concept of spinning reserve width and related results), spinning reserve constraints are also represented by (7). With some linearization techniques, the above SCUC problem can be transformed into a mixed integer linear programming problem (MILP) [15].

## III. VARIABLE REDUCTION APPROACH

The key in the proposed variable reduction approach is to quickly obtain an approximate UC solution, and then based on this solution, transform the original **f-UC** into a smaller scale SCUC problem to obtain the final UC solution. The approximate solution is obtained based on the detailed analysis of the single-period SCUC problem.

### A. Single-Period SCUC (s-UC)

The single-period SCUC problem is obtained by relaxing the time-coupled constraint (6) and $z \in X$ in (7) of the multi-period **f-UC**. The start-up cost is also neglected in the single-period SCUC problem, since it is related to the decision variables of consecutive time periods.

The single-period SCUC at period $t$ is defined as (8)-(13):

(**s-UC**):

$$\min_{z_t,p_t} \sum_i C_i(p_{i,t},z_{i,t}) \qquad (8)$$

$$\text{s.t. } (\lambda_{0,t}) \sum_i p_{i,t} - \sum_k d_{k,t} = 0 \qquad (9)$$

$$(\lambda_{l,t}^+) \sum_i \Gamma_{l,i}^U p_{i,t} - \sum_k \Gamma_{l,k}^D d_{k,t} \le F_l; \forall l \qquad (10)$$

$$(\lambda_{l,t}^-) - \sum_i \Gamma_{l,i}^U p_{i,t} + \sum_k \Gamma_{l,k}^D d_{k,t} \le F_l; \forall l \qquad (11)$$

$$z_{i,t}\underline{P} \le p_{i,t} \le z_{i,t}\overline{P}_i; \forall i \qquad (12)$$

$$z_t \in \{0,1\}^{I \times 1} \qquad (13)$$

Suppose the optimal UC solution of **s-UC** is $z'_t \in \{0,1\}^{I \times 1}$, then $z' = (z'_1, z'_2,...,z'_T) \in \{0,1\}^{I \times T}$ can be seen as a good starting point to obtain the final solution for the **f-UC** problem.

Though the scale of **s-UC** is much smaller than **f-UC**, there are two concerns: 1) **s-UC** is still a computationally demanding MILP problem. 2) $z'$ may not be a *feasible* solution to **f-UC**, and even if it is, the optimality of the solution is not guaranteed. The two concerns are addressed in subsection B and C. Afterward, subsection D gives the overall algorithm.

### B. Approximate Solution of s-UC by Linear Programming

#### 1) Lagrangian Relaxation of s-UC

For **s-UC** at period $t$, by introducing dual variable $\lambda_{0,t}$ to constraint (9) and $\lambda_{l,t}^+, \lambda_{l,t}^- \ge 0$ to constraints (10)-(11), the Lagrangian of the **s-UC** at period $t$ is as (14):

$$L_t(\lambda_t,p_t,z_t) = \sum_i C_i(p_{i,t},z_{i,t}) + \lambda_{0,t}(\sum_k d_{k,t} - \sum_i p_{i,t})$$
$$+ \sum_l \lambda_{l,t}^+(\sum_i \Gamma_{l,i}^U p_{i,t} - \sum_k \Gamma_{l,k}^D d_{k,t} - F_l) \qquad (14)$$
$$- \sum_l \lambda_{l,t}^-(\sum_i \Gamma_{l,i}^U p_{i,t} - \sum_k \Gamma_{l,k}^D d_{k,t} + F_l)$$

All dual variables are represented together by the vector $\lambda_t$. Then, the dual function of **s-UC** is the minimization of the Lagrangian (14) with respect to $p_t, z_t$, as (15)-(17):

$$\Phi_t(\lambda_t) = \min_{p_t,z_t} L_t(\lambda_t,p_t,z_t) \qquad (15)$$

$$\text{s.t. } z_{i,t}\underline{P}_i \le p_{i,t} \le z_{i,t}\overline{P}_i; \forall i \qquad (16)$$

$$z_t \in \{0,1\}^{I \times 1} \qquad (17)$$

Then, the dual problem (d-UC) of **s-UC** is the maximization of the dual function $\Phi_t(\lambda_t)$ with respect to $\lambda_t$, as (18)-(19).

(**d-UC**):

$$\max_{\lambda_t} \Phi_t(\lambda_t) \qquad (18)$$

$$\text{s.t. } \lambda_{l,t}^+, \lambda_{l,t}^- \ge 0; \forall l \qquad (19)$$

It has been proved that $\Phi(\lambda_t)$ is a concave function of $\lambda_t$, so **d-UC** is a convex optimization problem. The optimal solution to **d-UC** can provide us with useful information on the optimal solution to **s-UC**. The optimal solution to **d-UC** is usually obtained by using the iterative subgradient-based algorithms in literature, and only linear convergence rate is guaranteed. Zig-zagging and solution oscillation are often encountered when using the subgradient-based algorithms [7].

In this paper, we found that the optimal solution to **d-UC** can be obtained by solving a very simple linear program (LP) without using the iterative subgradient-based algorithms. Hence, the exact global optimal solution to **d-UC** can be obtained very efficiently. The transformation to LP is based on the separability of **d-UC** and the analytical expression of the solution to individual-unit subproblem.

*2) Individual-unit subproblem*

With given $\lambda_t$, (15)-(17) can be solved by solving $I$ individual-unit subproblems. (15)-(17) is reformulated as (20)-(22) for a better description.

$$\Phi_t(\lambda_t) = \min_{p_{i,t}, z_{i,t}} \sum_i L_{i,t}(\lambda_t, p_{i,t}, z_{i,t}) + \lambda_{0,t} \sum_k d_{k,t}$$
$$+ \sum_l (\lambda_{l,t}^+ g_{l,t}^+ + \lambda_{l,t}^- g_{l,t}^-) \quad (20)$$

$$\text{s.t. } z_{i,t}\underline{P}_i \le p_{i,t} \le z_{i,t}\overline{P}_i; \forall i \quad (21)$$

$$z_{i,t} \in \{0,1\}^{I \times 1} \quad (22)$$

where

$$L_{i,t}(\lambda_t, p_{i,t}, z_{i,t}) = C_i(p_{i,t}, z_{i,t}) + \beta_{i,t}(\lambda_t) p_{i,t} \quad (23)$$

$$\beta_{i,t}(\lambda_t) = -\lambda_{0,t} + \sum_l (\lambda_{l,t}^+ - \lambda_{l,t}^-) \Gamma_{l,i}^U \quad (24)$$

$$g_{l,t}^+ = -\sum_k \Gamma_{l,k}^D d_{k,t} - F_l \ , \ g_{l,t}^- = \sum_k \Gamma_{l,k}^D d_{k,t} - F_l \ ; \forall l \quad (25)$$

For each unit $i$ at period $t$, with given $\lambda_t$, the individual-unit subproblem (i-UC) is as (26)-(28),

(**i-UC**):

$$L_{i,t}(\lambda_t) \triangleq \min_{p_{i,t}, z_{i,t}} L_{i,t}(\lambda_t, p_{i,t}, z_{i,t}) \quad (26)$$

$$\text{s.t. } z_{i,t}\underline{P}_i \le p_{i,t} \le z_{i,t}\overline{P}_i \quad (27)$$

$$z_{i,t} \in \{0,1\} \quad (28)$$

*Theorem 1*: With given $\lambda_t$, the optimal objective value of **i-UC** ($L_{i,t}(\lambda_t)$ in (26)) is a piecewise linear concave function of $\beta_{i,t}(\lambda_t)$ (defined in (24)) as (29)-(30).

$$L_{i,t}(\lambda_t) = \min\{L_{i,t}^{C1}(\lambda_t), L_{i,t}^{C2}(\lambda_t), 0\} \quad (29)$$

$$\begin{cases} L_{i,t}^{C1}(\lambda_t) = [a_i + \beta_{i,t}(\lambda_t)]\overline{P}_i + b_i \\ L_{i,t}^{C2}(\lambda_t) = [a_i + \beta_{i,t}(\lambda_t)]\underline{P}_i + b_i \end{cases} \quad (30)$$

The optimal UC solution to **i-UC** is as (31)

$$z_{i,t}(\lambda_t) = \begin{cases} 0; & \text{if } L_{i,t}(\lambda_t) = 0 \\ 1; & \text{otherwise} \end{cases} \quad (31)$$

*Proof:* According to (26) and (28), we have:

$$L_{i,t}(\lambda_t) = \min\left\{\min_{p_{i,t}} L_{i,t}(\lambda_t, p_{i,t}, 0), \min_{p_{i,t}} L_{i,t}(\lambda_t, p_{i,t}, 1)\right\} \quad (32)$$

Based on (27), if $z_{i,t}$=0 then $p_{i,t} = 0$. Therefore based on (2) and (23) we know:

$$L_{i,t}(\lambda_t, p_{i,t}, 0) = L_{i,t}(\lambda_t, 0, 0) = 0 \quad (33)$$

For $L_{i,t}(\lambda_t, p_{i,t}, 1)$, based on (2) we have:

$$L_{i,t}(\lambda_t, p_{i,t}, 1) = C_i(p_{i,t}, 1) + \beta_{i,t}(\lambda_t) p_{i,t} = \left[a_i + \beta_{i,t}(\lambda_t)\right] p_{i,t} + b_i \quad (34)$$

Define $\hat{L}_{i,t}(\lambda_t)$ as (35).

$$\hat{L}_{i,t}(\lambda_t) \triangleq \min_{\underline{P}_i \le p_{i,t} \le \overline{P}_i} L_{i,t}(\lambda_t, p_{i,t}, 1) \quad (35)$$

For (34)-(35), two cases are possible:

*Case 1*: $\beta_{i,t}(\lambda_t) \le -a_i$. Then $L_{i,t}(\lambda_t, p_{i,t}, 1)$ is a monotonically decreasing function of $p_{i,t}$. Then the optimal solution to (35) is $p_{i,t}^* = \overline{P}_i$. And the optimal objective value is $\hat{L}_{i,t}(\lambda_t) = [a_i + \beta_{i,t}(\lambda_t)]\overline{P}_i + b_i$.

*Case 2*: $\beta_{i,t}(\lambda_t) > -a_i$. Then $L_{i,t}(\lambda_t, p_{i,t}, 1)$ is a monotonically increasing function of $p_{i,t}$. So the optimal solution to (35) is $p_{i,t}^* = \underline{P}_i$. And the optimal objective value is $\hat{L}_{i,t}(\lambda_t) = [a_i + \beta_{i,t}(\lambda_t)]\underline{P}_i + b_i$

Substitute the above optimal objective value and (33), (35) into (32), then (29)-(31) are obtained. Q.E.D.

Theorem 1 is based on the linear fuel cost function in (2). The related result with piecewise linear fuel cost function is presented in the Appendix. Theorem 1 provides the analytical expression of the solution to the **i-UC** and is important in transforming the **d-UC** into a linear program.

*3) Optimal solutions of d-UC by Linear Programming*

With the above results, the **d-UC** (18)-(19) is now a single-level optimization problem as (36)-(37).

$$\max_{\lambda} \Phi_t(\lambda_t) = \sum_i L_{i,t}(\lambda_t) + \lambda_{0,t} \sum_k d_{k,t} + \sum_l (\lambda_{l,t}^+ g_{l,t}^+ + \lambda_{l,t}^- g_{l,t}^-) \quad (36)$$

$$\text{s.t. } \lambda_{l,t}^+, \lambda_{l,t}^- \ge 0; \forall l \quad (37)$$

Where $L_{i,t}(\lambda_t)$ in (36) is a piecewise linear function of $\beta_{i,t}(\lambda_t)$ as (29)-(30) given in Theorem 1. Based on (29)-(30), the **d-UC** (36)-(37) can be equivalently transformed into a linear program by introducing a group of auxiliary variables $y_{i,t}$. The final formulation is given by (38)-(42).

(**d-UC-LP**)

$$\max_{\lambda, y_t} \Phi_t(\lambda_t) = \sum_i y_i + \lambda_{0,t} \sum_k d_{k,t} + \sum_l (\lambda_{l,t}^+ g_{l,t}^+ + \lambda_{l,t}^- g_{l,t}^-) \quad (38)$$

$$\text{s.t. } \lambda_{l,t}^+, \lambda_{l,t}^- \ge 0; \forall l \quad (39)$$

$$y_{i,t} \le 0; \forall i \quad (40)$$

$$y_{i,t} \le [a_i + \beta_{i,t}(\lambda_t)]\overline{P}_i + b_i; \forall i \quad (41)$$

$$y_{i,t} \le [a_i + \beta_{i,t}(\lambda_t)]\underline{P}_i + b_i; \forall i \quad (42)$$

By solving **d-UC-LP**, the optimal solution of **d-UC** is obtained. Suppose for period $t$, the optimal solution of **d-UC-LP** is ($\tilde{\lambda}_t$, $\tilde{y}_t$). If $\tilde{y}_{i,t}$=0, then according to (31) in Theorem 1, the corresponding UC solution to (26)-(28) (for $\lambda_t = \tilde{\lambda}_t$) is $\tilde{z}_{i,t}$=0. Otherwise, $\tilde{z}_{i,t}$=1. This UC solution $\tilde{z}_t \in \{0,1\}^{I \times 1}$ is then used in variable reduction.

*C. Variable Reduction*

The UC solution $\tilde{z}$ may not be a feasible solution to the **f-UC** and, even if it is, it may deviate from the optimal solution of **f-UC**. However, we found that high quality UC solutions to **f-UC** can be obtained based on $\tilde{z}$. The basic idea is to solve a small scale SCUC by fixing some binary variables in the **f-UC** based on the information of the obtained $\tilde{\lambda}$, $\tilde{y}$ and $\tilde{z}$.

For unit $i$ at period $t$, based on (29)-(30), the figure of $L_{i,t}(\lambda_t)$ over $\beta_{i,t}(\lambda_t)$ is depicted as Fig. 2.

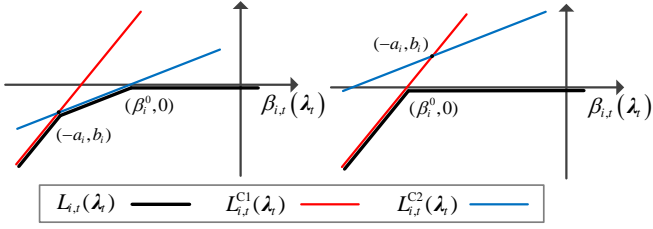

Fig. 2  $L_{i,t}(\lambda_t)$ over $\beta_{i,t}(\lambda_t)$ when $b_i<0$ (left) and $b_i \geq 0$ (right).

In Fig. 2, the value of intersection point of $\beta_i^0$ can be calculated as (43).

$$\beta_i^0 = \max\{-a_i - b_i \,/\, \overline{P}, -a_i - b_i \,/\, \underline{P}\} \qquad (43)$$

According to Fig. 2 and Theorem 1: 1) if $\beta_{i,t}(\lambda_t) > \beta_i^0$, then $L_{i,t}(\lambda_t)$ =0, and hence $z_{i,t}(\lambda_t)$ =0. 2) if $\beta_{i,t}(\lambda_t) \leq \beta_i^0$, then $z_{i,t}(\lambda_t)$ =1.

It is concluded from the above analysis that whether or not $\beta_{i,t}(\lambda_t)$ is less than $\beta_i^0$ is an indicator of whether the UC decision is on or off.

We then choose the distance from $\beta_{i,t}(\lambda_t)$ to $\beta_i^0$, i.e. $|\beta_{i,t}(\tilde{\lambda}_t) - \beta_i^0|$, as a criterion to fix the UC variables: for those UC variables with large $|\beta_{i,t}(\tilde{\lambda}_t) - \beta_i^0|$ should be fixed to the value of $\tilde{z}_{i,t}$, while the others are still binary variables to be determined.

### D. Overall Algorithm

The overall algorithm is as follows:

**(Algorithm 1)**

***Initialization.*** Set "Pre-Determined Ratio" (PDR).

***Step 1*** For each period $t$, solve **d-UC-LP** (38)-(42), obtain the optimal solution $\tilde{\lambda}_t$, $\tilde{y}_t$.
For each $i,t$, If $\tilde{y}_{i,t}=0$, then set $\tilde{z}_{i,t}$ =0; else, $\tilde{z}_{i,t}$ =1.

***Step 2*** For each unit $i$, calculate $\beta_i^0$ based on (43). For each $i$ and period $t$, calculate $\beta_{i,t}(\tilde{\lambda})$ based on (24), and then calculate $|\beta_{i,t}(\tilde{\lambda}) - \beta_{i,t}^0|$.

***Step 3*** Sort $|\beta_{i,t}(\tilde{\lambda}_t) - \beta_{i,t}^0|$ of all $i,t$ in descending order:
$|\beta_{i_1,t_1}(\tilde{\lambda}_t) - \beta_{i_1,t_1}^0| > ... > |\beta_{i_T,t_T}(\tilde{\lambda}_t) - \beta_{i_T,t_T}^0|$

***Step 4*** Suppose $I \times T \times$PDR$=N$, and $z$ represents the UC variable of **f-UC**. For each $n=1,...,N$, replace variable $z_{i_n,t_n}$ with the 0/1 constant $\tilde{z}_{i_n,t_n}$.

***Step 5*** Solve the new **f-UC** with the reduced variables by MIP solver. If no feasible solutions can be found, then lower the PDR, and go to ***Step 4***. Otherwise, stop and return the optimal solution.

In the initialization step, "Pre-Determined Ratio" (PDR) is a parameter set by users, which represents the percentage of the UC variables to be pre-determined according to $\tilde{z}$. When the PDR increases, more UC variables in **f-UC** are replaced with 0/1 constants, and the problem scale becomes smaller. However, when the PDR is too large, the pre-determined UC decisions may violate some time-coupled constraints, so the problem can be infeasible. In this case, reducing the PDR can help to find the feasible solution.

The variable reduction technique proposed in this paper can be used together with many other techniques for large-scale SCUC problems like the technique of fast elimination of the redundant security constraints [21], locally ideal formulation proposed in [16], and so on.

## IV. NUMERICAL RESULTS

All numerical tests are performed with MATLAB R2015b, YALMIP toolbox [23] and GUROBI 8.0 on an Intel Core(TM) i7-3770 CPU @ 3.40GHz PC with 20GB RAM.

The numerical tests are conducted in two cases: 1) the IEEE 118-bus system with N-1 security constraints [24]. 2) the modified French 6468-bus system with base case security constraints [25]. Necessary information on these two systems are in TABLE I. Fuel cost functions of both cases are piecewise-linear with three segments.

TABLE I Basic Information of The Test Systems

| System | 118-bus | 6468-bus |
|---|---|---|
| # of Units | 54 | 399 |
| # of Periods | 24 | 24 |
| # of Transmission Lines | 179 | 9000 |
| N-1 Constraints | Included | Not Included |
| # of Transmission/Security Constraints | 1366128 | 432000 |
| Generation Capacity (MW) | 13374 | 115596 |
| Peak Load (MW) | 6000 | 52535 |
| MIP Gap for Solver | 0.01% | 0.05% |

The redundant security constraints are first eliminated by the method proposed in [21]. In the 118-bus system, 823382 out of 1366128 (60.27%) security constraints are removed from the model. In the 6468-bus system, 374029 out of 432000 (86.58%) security constraints are removed. The Elimination of the redundant security constraints is very helpful in reducing the RAM usage and improving the computational efficiency.

The rest of this section is organized as follows: In subsection A, the solution quality and computational efficiency of **d-UC-LP** is tested and compared with **s-UC** and **f-UC** (corresponds to section III.A and III.B.). In subsection B, the effectiveness of the criterion for variable selection and reduction is tested (corresponds to section III.C.). In subsection C, the overall variable reduction technique is tested (corresponds to the whole section III).

### A. d-UC-LP: Solution Quality and Computational Efficiency

In this subsection, the test cases are solved by 1) **f-UC** (1)-(7) to obtain the benchmark UC solution, 2) **s-UC** (8)-(13)(for 24 time periods) to see if its solution is similar to the solution of **f-UC**, and 3) the **d-UC-LP** (38)-(42) (for 24 time periods) It is again noted that the solutions of **s-UC** and **d-UC-LP** are only the approximate UC solutions or even infeasible to **f-UC**, MIP solvers must be launched afterward to obtain the final UC solution, as stated in section III.D.

The results for 118-bus and 6468-bus systems are presented in TABLE II and TABLE III, respectively. The start-up cost of **f-UC** is subtracted from the optimal objective value, and 24 optimal objective values of **s-UC** and **d-UC-LP** (each corresponds to one period) are added together for comparison. "# of committed units" represents the total number of committed units along the 24 time periods. "# of Different UC" represents the number of UC decisions that are different from **f-UC**.

TABLE II Results of **f-UC**, **s-UC**, **d-UC-LP** on 118-bus System

| Methods | **f-UC** | **s-UC** | **d-UC-LP** |
|---|---|---|---|
| Objective Value (E6$) | 1.8085 | 1.8080 | 1.7980 |
| # of Committed Units | 655 | 644 | 625 |
| # of Different UC | 0 | 21 | 72(59 to s-UC) |

TABLE III Results of **f-UC**, **s-UC**, **d-UC-LP** on 6468-bus System

| Methods | **f-UC** | **s-UC** | **d-UC-LP** |
|---|---|---|---|
| Objective Value (E7$) | 1.0648 | 1.0610 | 1.0597 |
| # of Committed Units | 1128 | 1078 | 1091 |
| # of Different UC | 0 | 284 | 323(165 to s-UC) |

TABLE IV $\Delta\beta_{i,t}$ of elements with $\tilde{z}_{i,t} = z'_{i,t}$ and $\tilde{z}_{i,t} \neq z'_{i,t}$

| Elements | 118-bus | | 6468-bus | |
|---|---|---|---|---|
| | # | Avg. $\Delta\beta_{i,t}$ | # | Avg. $\Delta\beta_{i,t}$ |
| all $i,t$ | 1296 | 5.3154 | 9576 | 14.6056 |
| $\tilde{z}_{i,t} = z'_{i,t}$ | 1237 | 5.5630 | 9411 | 14.8614 |
| $\tilde{z}_{i,t} \neq z'_{i,t}$ | 59 | 0.1229 | 165 | 0.012 |

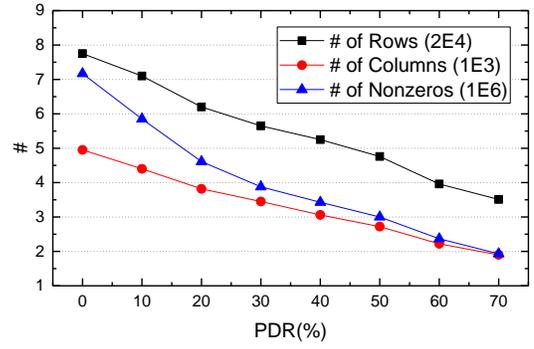

Fig. 3  # of rows, columns and nonzero elements after presolve with different PDR settings in 118-bus case

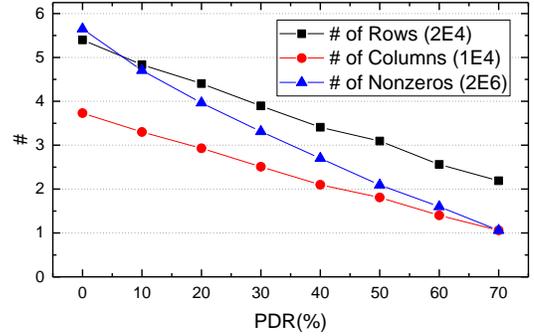

Fig. 4 # of rows, columns and nonzero elements after presolve with different PDR settings in 6468-bus case

It is seen from TABLE II and TABLE III that:

1) The objective values satisfy **f-UC** > **s-UC** > **d-UC-LP** for both cases. For the first inequality, the gap is 0.03% and 0.36%, respectively for the two cases. This is because **f-UC** is constrained by time-coupled constraints while **s-UC** is not, so the feasible region of **f-UC** is a subset of that of the **s-UC**, and the objective value is thus larger. For the second inequality, the gaps of the two cases are 0.55% and 0.12%, respectively.

2) The UC decisions obtained by all methods are close to each other, according to both the number of committed units and the number of different UC decisions. It is verified that the solution of **s-UC** is similar to the solution of **f-UC** (as stated in section III.A.). In the 118-bus and 6468-bus cases, only 1.62% (21 out of 1296) and 2.97% (284 out of 9576) UC decisions of the **s-UC** are different from the **f-UC**.

The difference between the UC solutions of **d-UC-LP** and **s-UC** is also small, as only 4.55% (59 out of 1296) and 1.72% (165 out of 9576) of the UC decisions are different. The results suggest that the idea of solving **d-UC-LP** to obtain the approximate solution of **s-UC** (stated in section III.B.) is also effective.

Regarding the computational time, **s-UC** uses 177.4s and 835.1s, respectively for the two cases, to obtain the optimal solution. While **d-UC-LP** reduces the time to 45.2 and 45.7 seconds. This is because **d-UC-LP** is a linear programming problem, while **s-UC** is still an MILP problem. Therefore, **d-UC-LP** is suitable for a preprocessing procedure regarding the computational effort.

### B. The Criterion for Variable Selection and Reduction

In this subsection, we test the criterion for variable selection and reduction proposed in section III.C.

Denote the optimal UC solutions of **d-UC-LP** and **s-UC** as $\tilde{z}$ and $z'$, respectively. And denote:

$$\Delta\beta_{i,t} = |\beta_{i,t}(\tilde{\lambda}_t) - \beta^0_{i,t}| \quad (\forall i,t) \tag{44}$$

where $\tilde{\lambda}_t$ is the optimal solution of **d-UC-LP**. $\beta_{i,t}(\tilde{\lambda}_t)$ and $\beta^0_{i,t}$ are obtained in Step 2 in Algorithm 1.

The idea of the criterion is to select the UC variables with large $\Delta\beta_{i,t}$, and then fix them to the values of $\tilde{z}_{i,t}$, since we believe the solution of $z'_{i,t}$ with large $\Delta\beta_{i,t}$ is more likely to be the same to $\tilde{z}_{i,t}$. Therefore, we compare the average $\Delta\beta_{i,t}$ of 1) elements with $\tilde{z}_{i,t} = z'_{i,t}$ and 2) the elements with $\tilde{z}_{i,t} \neq z'_{i,t}$. The results are in TABLE IV.

It is seen from TABLE IV that: For the elements with $\tilde{z}_{i,t} = z'_{i,t}$, the average $\Delta\beta_{i,t}$ is much larger than those with $\tilde{z}_{i,t} \neq z'_{i,t}$. The results suggest that $\Delta\beta_{i,t}$ is a very good criterion that can be used to identify the UC decisions of **d-UC-LP** that are likely to be the same with **s-UC**.

### C. Overall Variable Reduction approach.

In this subsection, we test the overall variable reduction technique stated in Algorithm 1. Algorithm 1 is executed with PDR set from 0% to 100% with 10% as the step size (therefore, Algorithm 1 is executed for 11 times for each case).

When PDR is 0%, no variables are replaced in **f-UC**, the problem is the same with the full-scale **f-UC**, and the results are seen as the benchmark**.** When PDR increase to 80%, no feasible solutions can be found for both cases. As analyzed in section III.C, the solution of **s-UC** may not be the feasible solutions of **f-UC** due to the relaxation of time-coupled constraints by **s-UC**. In the rest of this subsection, we only present the results with PDR=0%~70%.

The number of integer (binary UC) variables is a critical factor that affects the potential scale of the branch-and-cut tree

[16], and thus affects the overall computational burden. In 118-bus case, the number of the integer variables is reduced from 1296 to 389 as PDR increases from 0% to 70%, and 9676 to 2873 in the 6468-bus case.

Before solving the problem, the MIP solver "presolves" the problem to "make a model smaller and easier to solve" [26] by removing unnecessary rows and columns (variables) in the co-efficient matrix. With the reduced integer variables, the number of rows, columns and nonzero elements after presolving are presented in Fig. 3 and Fig. 4. It is seen from the results that 1) in both cases, the numbers of rows, columns, and nonzeros monotonically decrease as PDR gets larger. 2) The number of columns equals the number of the variables, as the integer variables are reduced, the columns are also reduced. 3) The number of rows decreases to 45.14% and 40.43% of the full-scale model, respectively in the two cases, when PDR increases from 0% to 70%. As more variables are replaced by 0/1 constants, some constraints may become redundant and are removed by the solver. 4) The number of nonzero elements in the coefficient matrix is also a key factor that can affect the computational burden [3], which decrease to 26.86% and 18.68% of the full-scale model, respectively in two cases. As rows and columns are reduced, the nonzero elements within the coefficient matrix are naturally reduced.

With different PDR settings, the **f-UC** with reduced integer variables are now solved by the MIP solver, and the results are presented in TABLE V and TABLE VI. "CPU Time" is the solver time returned by the solver. "Obj. Value" is the optimal objective value. "Root Relax." is the time spent by root relaxation. At each node of the branch-and-cut tree, an LP relaxation problem is solved, with a similar structure to the relaxed problem of other nodes, so the time spent by root relaxation can help to evaluate the computational burden of these LP relaxation problems. "Nodes Expl." is the number of nodes explored by the solver in the branch-and-cut tree.

It is seen from TABLE V and TABLE VI that:

1) In both cases, the computational time spent by the solver generally decreases as PDR increases. Moreover, when PDR is 70%, it takes only 6.13% and 29.47% of the time spent to solve the full-scale **f-UC** (when PDR=0%).

2) In 6468-bus case, the objective values with different PDR settings are the same, which means that when up to 70% of the UC variables are fixed, the optimality of the solution is not affected. In 118-bus case, the objective value of PDR=60% and PDR=70% are 0.01% and 0.15% larger than the rest of the PDR settings, so the optimality of the solution are affected, but to a very limited extent.

3) The time spent by root relaxation generally decreases as PDR increases, which indicates that it is probably easier to solve the node relaxation problems when more variables are reduced. This can also be explained by the results presented in Fig. 3 and Fig. 4, since the node relaxation problems have the similar coefficient matrices that are derived from the matrix whose scale is presented in the two figures.

TABLE V f-UC with Variable Reduction of 118-bus System

| PDR | CPU Time (s) | Obj. Value (E6$) | Root Relax. (s) | Nodes Expl. |
|---|---|---|---|---|
| 0 | 1987.6 | 1.8092 | 8.49 | 6290 |
| 10% | 1588.5 | 1.8092 | 6.84 | 4998 |
| 20% | 547.4 | 1.8092 | 5.27 | 3214 |
| 30% | 811.9 | 1.8092 | 6.29 | 5254 |
| 40% | 757.7 | 1.8092 | 5.67 | 5389 |
| 50% | 478.0 | 1.8092 | 4.47 | 5119 |
| 60% | 200.3 | 1.8094 | 3.65 | 1461 |
| 70% | 121.9 | 1.8120 | 2.27 | 455 |

TABLE VI f-UC with Variable Reduction of 6468-bus System

| PDR | CPU Time (s) | Obj. Value (E7$) | Root Relax. (s) | Nodes Expl. |
|---|---|---|---|---|
| 0 | 4255.2 | 1.0668 | 94.11 | 2371 |
| 10% | 3632.1 | 1.0668 | 133.91 | 2456 |
| 20% | 2368.8 | 1.0668 | 69.36 | 1937 |
| 30% | 2240.8 | 1.0668 | 63.48 | 1587 |
| 40% | 1956.0 | 1.0668 | 92.84 | 2084 |
| 50% | 1864.4 | 1.0668 | 51.93 | 2662 |
| 60% | 1430.0 | 1.0668 | 41.62 | 2137 |
| 70% | 1254.1 | 1.0668 | 38.49 | 1911 |

4) In 118-bus case, the nodes explored when optimality of the solution is affected (PDR=60% and 70% in 118-bus case) are much fewer than those settings where the optimality of the solution is not affected. While in 6468-case, the numbers of nodes explored are close to each other with different PDR settings. However, with the analysis presented in 3), the overall computational time still decreases as the solver may spend less time exploring each node.

## V. CONCLUSIONS

In this paper, a variable reduction method is proposed and applied in large-scale SCUC problem. The solution of the single-period **s-UC** problem is close to the full-scale multi-period **f-UC** problem, and can be seen as a good starting point to obtain the final UC solution for **f-UC**. The approximate solution of **s-UC** as well as its optimal dual solution can be directly obtained by solving a linear program **d-UC-LP**; the optimality gap is small and the computational time is negligible compared to directly solving **s-UC**. A criterion of selecting binary variables is established based on the solution of **d-UC-LP**. The selected binary variables are then fixed to the value of the approximate solution of **s-UC**. With reduced binary variables, the **f-UC** has much smaller problem scale, and the computational efficiency improves significantly, without sacrificing too much of its optimality.

## Appendix

### Theorem 1 with Piecewise Linear Fuel Cost Functions

Suppose for any unit $i$, the fuel cost function is as follows:

If $z_{i,t} = 0$, then $C_i(p_{i,t},0) = 0$

If $z_{i,t} = 1$, then

$$C_i(p_{i,t},1) = \begin{cases} a_{i,1}p_{i,t} + b_{i,1}; & \text{if } \underline{P} \le p_{i,t} \le P_{i,1} \\ a_{i,2}p_{i,t} + b_{i,2}; & \text{if } P_{i,1} \le p_{i,t} \le P_{i,2} \\ \dots \\ a_{i,M_i}p_{i,t} + b_{i,M_i}; & \text{if } P_{i,M_i-1} \le p_{i,t} \le \overline{P}_i \end{cases}$$

where $\underline{P}_i = P_{i,0} < P_{i,1} < \dots < P_{i,M_i-1} < P_{i,M_i} = \overline{P}_i$, and it is assumed that $a_{i,1} \le a_{i,2} \le \dots \le a_{i,M_i}$ ($C_i(\cdot)$ is convex w.r.t $p_{i,t}$).

Then the conclusions in Theorem 1 are as follows:

*Theorem 2*: the optimal objective value of **i-UC** $L_{i,t}(\lambda_t)$ is a piecewise linear concave function of $\beta_{i,t}(\lambda_t)$ as (45)-(46), shown in Fig. 5.

$$L_{i,t}(\lambda_t) = \min\{\hat{L}_{i,t}(\lambda_t), 0\} \tag{45}$$

$$\hat{L}_{i,t}(\lambda_t) = \begin{cases} [a_{i,M_i} + \beta_{i,t}(\lambda_t)]\overline{P}_i + b_{i,M_i}; \text{if } \beta_{i,t}(\lambda_t) \le -a_{i,M_i} \\ [a_{i,M_i-1} + \beta_{i,t}(\lambda_t)]P_{i,M_i-1} + b_{i,M_i-1}; \\ \qquad\qquad \text{if } -a_{i,M_i} < \beta_{i,t}(\lambda_t) \le -a_{i,M_i-1} \\ \dots \\ [a_{i,1} + \beta_{i,t}(\lambda_t)]P_{i,1} + b_{i,1}; \text{ if } -a_{i,2} < \beta_{i,t}(\lambda_t) \le -a_{i,1} \\ [a_{i,1} + \beta_{i,t}(\lambda_t)]\underline{P}_i + b_{i,1}; \text{ if } -a_{i,1} < \beta_{i,t}(\lambda_t) \end{cases} \tag{46}$$

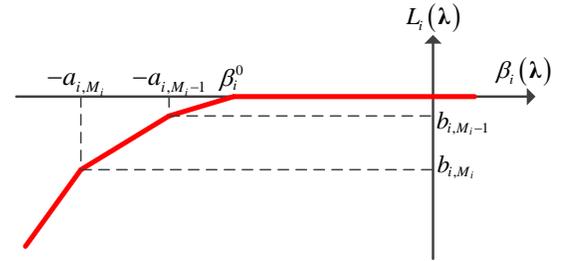

Fig. 5 Optimal objective value of **i-UC** with piecewise linear fuel cost functions. The figure shows the case when $b_{i,M_i-1} < 0 < b_{i,M_i-2}$.

The proof is very similar to Theorem 1, and is omitted due to length limit.